\documentclass{article}
\usepackage{amsfonts,amssymb,graphicx}
\newtheorem{theorem}{Theorem}[section]

\newtheorem{lemma}[theorem]{Lemma}

\newtheorem{corollary}[theorem]{Corollary}
\newtheorem{question}{Question}
\newcommand{\bproof}{\paragraph {\bf Proof. }}
\newcommand{\qed}{\square}

\begin{document}
\title{Boundary slopes of immersed surfaces in 3-manifolds}
\author {Joel Hass, \and J.Hyam Rubinstein, \and Shicheng Wang}

\maketitle
\begin{abstract} 
This paper presents some finiteness results for the number of
boundary slopes of immersed proper $\pi_1$-injective surfaces of given genus $g$ in a compact 3-manifold with torus boundary.
In the case of hyperbolic 3-manifolds we obtain uniform quadratic bounds in
$g$, independent of the 3-manifold.
\end{abstract}


\section{Introduction.}

An immersed, proper, $\pi_1$-injective
surface in a compact 3-manifold $M$ with non-empty boundary is {\it essential} if it cannot be properly homotoped into $\partial M$.
Let $c$ be a homotopically non-trivial simple loop in $\partial M$.
If there is a proper immersion of an essential surface $F$
into $M$ such that each component of $\partial F$
is homotopic to a multiple of $c$,
we call $c$ a {\it boundary slope} of $M$.
The first question we look at is a problem of P. Shalen, told to us by M. Baker:

\begin{question} \label{q1}
Does the set of essential surfaces with bounded genus in a simple
knot complement give rise to at most finitely many boundary slopes?
\end{question}

Baker has given examples to show that if the bounded genus assumption is
dropped then infinitely many boundary slopes can be realized \cite{Ba},
and Oertel has found examples of manifolds in which every slope is
realized by the boundary of an immersed essential surface \cite{Oe}, see also
Maher \cite{Maher}.
On the other hand,
Hatcher \cite{Ht} has shown that there are only finitely
many boundary slopes for embedded essential surfaces,
without a genus restriction. 

We answer Question~\ref{q1} in Sections 4 and 5, proving a stronger result.
Minimal surface theory is used (Theorem~\ref{4.3} and Theorem~\ref{4.4})
to derive a bound which is a quadratic function of $g$, independent of $M$,
in case the interior of $M$ has a complete hyperbolic metric of finite volume.
We also find a closely related upper bound $n(g_1,g_2)$ for the number of
intersections of $\alpha_1$ and $\alpha_2$, where $\alpha_i$ is a boundary
slope of an immersed surface of genus $g_i$, $i=1,2$.
With an additional combinatorial argument, a positive answer of the question for 3-manifolds with non-trivial
Jaco-Shalen-Johannson decomposition is given in Section 5. 

When $g=0$ or $1$, and the surface is an embedded punctured sphere or torus, there are many known results on the above questions, some sharp. These are based on highly developed combinatorial methods in knot theory and the theory of representations of knot groups.
See the survey papers \cite{Go}, \cite{Lu} and \cite{Sh}.
In the case where the surfaces are immersed punctured spheres or tori,
the Gromov-Thurston $2\pi$-lemma can be used to give bounds \cite{BH}.
The use of minimal surface theory to give uniform bounds for the number
of boundary slopes of $\pi_1$-injective immersed surfaces of bounded genus is natural,
and does not seem to have appeared before in this context,
though it is inspired by the work of Thurston, Uhlenbeck and Meeks-Yau.

Another question we investigate was raised by J. Luecke \cite{Lu},\cite{Lu1}.
Let $K$ be a simple knot in $S^3$ and $(K,\lambda )$ be the closed 3-manifold obtained by surgery on $K$ along a slope $\lambda$.
Let $c(K, g, \lambda)$ be the least upper bound for
the geometric intersection numbers of the
core of the surgery solid torus and the homotopy class of any
closed essential surface of genus $g$ in $(K,\lambda)$. 
So any closed essential surface of genus $g$ in $(K,\lambda)$ can be homotoped
to intersect the core of the surgery in at most $c(K, g, \lambda)$ points.

\begin{question} \label{q2}
Is there a universal upper bound $c(g)$ for all $c(K, g, \lambda)$,
independent of the choices of the knot $K$ and slope $\lambda$?
\end{question}

Luecke found such a bound for non-integral surgeries on $K$. Closely related
to this question is the study of a bound
$n(K, g,\lambda)$ for the number of boundary components of a surface $F$
of genus $g$ in a knot complement in $S^3$ with all components of $\partial F$ having the
same slope $\lambda$.
We study the relationship between $c(K,g,\lambda)$ and $n(K,g,\lambda)$ and get a
partial answer to Question~\ref{q2}.

To describe our bounds on the number of slopes, we need some terminology. 
Let $E^2$ be the Euclidean plane and let ${\cal G}_d$, $d \ge 1$ be the set of all lattices on $E^2$ satisfying:

\begin{enumerate}
\item any two vertices of the lattice have distance at least $d$, 

\item the area of each parallelogram which is a
fundamental domain for the lattice is at least $d^2\sqrt 3$,

\item the origin of $E^2$ is a vertex.

\end{enumerate}

Call a non-zero vertex of a lattice {\it primitive} if it is not a positive integral
multiple of any other vertex. Let $D(R)$ be the disc of radius $R$ and let $$
N(g,d)= \frac 12 \max_{\{\Gamma_d \in {\cal G}_d \}} \{ \mbox{the number of primitive vertices of } \Gamma_d \mbox{ in } D(2g\pi) \} $$
In our applications the constant $d$ 
will depend on the geometry of the cusp of a hyperbolic 3-manifold, but will always be at least one.
Note that $N$ is non-decreasing as a
function of $g$ and non-increasing as a function of $d$.
$N$ can also be viewed as a function of the single variable $\frac gd$.
When $d=1$, we write $ N(g)$ for $N(g,1)$, and note that $N(g,d) \le N(g)$.
For any given pair $(g,d)$, the value of $N(g,d)$ can be computed. We will show in Section 3 that
$N(1,1) = 24$, $N(2,1)=92$, $N(10,1)=2186$ and $N(2, 1.15)=69$. We establish in Theorem~\ref{3.5}
that $N(g)$ is bounded by a quadratic function of $g$. More precisely, we show:
$$
\displaystyle \lim_{g \to \infty}\frac
{N(g,d)} {4\sqrt {3} (g+0.5)^2\pi}\le 1. $$

In a hyperbolic 3-manifold with boundary a torus we show that: 

\begin{enumerate}
\item
The number of boundary slopes of
essential immersed surfaces of
genus smaller or equal to $g$ is at most $N(g,d)$ for some $d \ge 1$ (Theorem~\ref{4.1}), and
so this number grows at most quadratically with $g$. It follows that
if $M$ contains no closed
essential surfaces of genus at most $g$, then at most $N(g)+1$ surgeries on
the cusp give closed 3-manifolds
containing closed $\pi_1$-injective
surfaces with genus at most $g$.

\item There is a quadratic bound $n(g_1,g_2)$ for the intersection of
two slopes
$\alpha_1$ and $\alpha_2$, where $\alpha_i$
is the boundary slope of an immersed surface of genus $g_i>0$, $i=1,2$.
With some specified exceptions the bound is $11.8g_1g_2$;
for details see Theorem~\ref{4.5}.
C. Gordon has informed us that a combinatorial method developed by Gordon-Litherland establishes a quadratic bound for the special case where the surfaces are embedded.

\item $c(K, g, \lambda) \le n(K, g, \lambda)+1$ for any $K, g, \lambda$ and $n(K, g,\lambda) \le 2g-2$,
except for 92 possible exceptional slopes $\lambda$ (see Theorems~\ref{4.6} and~\ref{4.7}). C. Gordon has informed us that
in a hyperbolic knot complement in $S^3$, an embedded torus with an arbitrarily large number of punctures can be found. 
\end{enumerate}

For a general Haken manifold, such finiteness results also exist, but the bounds are not explicit, and depend on the manifold. 

The second section contains preliminary results.
In the third section we discuss the computation of $N(g,d)$.
We prove the main results in Section 4 and Section 5. 

All surfaces and 3-manifolds considered in this paper are assumed to be connected and orientable.

\section{Preliminaries.}

A map of a surface into a 3-manifold $f:F \to M$ is {\it $\pi_1$-injective} if the induced map
on the fundamental groups $f_* : \pi_1 (F) \to \pi_1 (M) $ is an injective homomorphism.
Given a subsurface $A \subset \partial M$ and a map $f: (F, \partial F) \to (M, A)$, we say that the surface is {\it $\pi_1$-injective relative to $A$} if $f_* : \pi_1 (F, \partial F) \to \pi_1 (M, A) $ is an injective homomorphism.
This means that any proper arc in $(F, \partial F)$
that has image which is homotopic to $A$ in $M$ (rel boundary) is homotopic to $\partial F$ in $F$ (rel boundary). 

An immersed $\pi_1$-injective surface which is not properly homotopic to the boundary of $M$ is
an {\it essential} immersed surface.
Any mapping of a surface into a 3-manifold is homotopic to an immersion in its interior,
by the classical construction of canceling
interior branch points or pushing them to the boundary of the 3-manifold.
It is not possible in general to perturb away boundary singularities.
For example, a figure eight on the plane in $R^3$ does not bound an immersed disk in the upper half-space.
In this paper we will consider immersions with no boundary singularities. 

We begin by examining the relationship between a surface which is injective on $\pi_1$ and one which is also injective on relative $\pi_1$. 

\begin{lemma} \label{lemma2.1}
Let $M$ be a compact irreducible 3-manifold with boundary,
and let $T$ be a torus boundary component of $M$.
Let $F$ be a $\pi_1$-injective surface with $\partial F \subset T$.
Then either $F$ is a boundary parallel annulus or $F$ is also injective on relative $\pi_1$.
\end{lemma}

\bproof Suppose that $F$ is a $\pi_1$-injective surface with
$\partial F \subset T$ and that $F$ is not injective on relative
$\pi_1$. Let $\alpha$ be an arc on $F$, not boundary parallel,
which is homotopic (rel boundary) into $T$. If $\alpha$ connects
two distinct boundary components $\beta$ and $\gamma$ of $\partial F$
then it follows that $\beta$ and $\alpha \gamma \alpha^{-1}$ are both
homotopic in $\pi_1 (M)$ into $T$, and thus they commute in $\pi_1 (M)$.
Since we assumed that $F$ is a $\pi_1$-injective surface,
it follows that $\beta$ and $\alpha \gamma \alpha^{-1}$ commute in $\pi_1(F)$.
But the elements of the fundamental group represented by two distinct boundary
components of a surface can commute only if the surface is an annulus,
and in this case boundary compressibility and irreducibility imply that
the annulus is boundary parallel.
If $\alpha$ connects a boundary component $\beta$ to itself,
then $\beta$ and $\alpha \beta \alpha^{-1}$ are both homotopic in $\pi_1 (M)$
into $T$, and thus commute in $\pi_1 (M)$.
Arguing as before, we conclude that a boundary component of $F$
commutes in $\pi_1(F)$ with a non-trivial conjugate of itself,
which is impossible. The lemma follows. $~~~\qed$

Let $M$ be a complete hyperbolic 3-manifold with finite volume and $\partial M$
a union of horotori. The boundary of $M$ has a flat
Riemannian metric induced from the hyperbolic metric. The {\it cusp length} of $M$, CL($M$), is the
supremum, over all choices of horotorus boundary, of the length of the shortest
Euclidean geodesic on a boundary torus of $M$. Adams \cite{Ad2} has made an extensive study of the cusp length. 
Adams showed that the figure eight knot complement plays a special role.
It alone can have a cusp length equal to one.

\begin{lemma} \label{lemma2.2}
The cusp length of $M$ satisfies CL($M )\ge 1$ for any hyperbolic
3-manifold with torus boundary components and
CL($M )\ge 1.15$ for any $M$ other than the complement of the figure eight knot.
Moreover the area of a maximal cusp is at least $CL(M)^2\sqrt 3$,
and the minimal area among all maximal cusps is at least $3.35$.
\end{lemma}

\bproof A horotorus cutting off a
cusp can be pushed into $M$ until it touches either itself or another cusp.
Hyperbolic geometry shows that the distance along the boundary torus between closest points of tangency of cusps is at
least one.
Adams analyzed the possible configurations with small cusp lengths,
and deduced that cusp length smaller than 1.15094 is only possible in the cases given above \cite{Ad2}.
The cusp area bound CL$(M)\sqrt 3$ is also due to Adams (Theorem 2, \cite{Ad1}),
and the bound 3.35 was recently found by Cao and Meyerhoff \cite {CM} $~~~\qed$
Now we state some facts about surfaces in Seifert manifolds.
First, it is known that each $\pi_1$-injective
surface in a Seifert manifold can be homotoped to be either vertical or horizontal \cite{Ha}.
Now suppose $p:M\to F$ is an oriented Seifert manifold, where $F$ has
genus $g$ and $h>0$ boundary components and $M$ has $k$ singular fibers. $M$ is {\it framed} if 

\begin{enumerate}
\item a section $S= F-\cup \mbox{int} (D_i)$ of $M-\cup \mbox{int} (N_i)$ is chosen and $\partial S$ is
oriented, where the $N_i$ are fibered regular neighborhoods of the singular fibers;

\item each torus boundary component $T_i$ of $M$ is equipped with a
framing $T_i (\mu_i ,\lambda_i )$, where $\mu_i$ is an oriented boundary
component of $F$ and $\lambda_i$ is an oriented fiber; 

\item each torus $T_i$ is given the orientation induced from $M$. 
\end{enumerate}

Once the section $S$ is chosen the data
$(g;h; \alpha_1,\beta_1;\dots ;\alpha_k,\beta_k)$ specifies the Seifert fiber space $M$.
Suppose that $T_i(\mu_i,\lambda_i)$,
$i=1,\dots ,n$ are the boundary components of the framed Seifert manifold $M$. 

\begin{lemma} \label{lemma2.3}
Let $F$ be an essential horizontal orientable immersed surface in an
orientable Seifert manifold $M$ and let
$\{c_{i,j}, j=1,\dots ,k_i\}$ be its boundary components (with induced orientation) on
$T_i$, $c_{i,j}=(u_{i,j},v_{i,j})$, $i=1,\dots ,h$.
Let $u$ be the geometric intersection number of $F$ and a regular fiber and let
$h_i$ be the number of boundary components of $F$ on $T_i$.

Then 
\begin{equation}
\sum_{j=1}^{h_i} u_{i,j}=u\ne 0 \label{eq2.1}
\end{equation} 

\begin{equation}
u\sum_{i=1}^k\frac {\beta_i}{\alpha_i}+ \sum_{i=1}^h\sum_{j=1}^{h_i}\frac {v_{i,j}} {u_{i,j}}=0
\label{eq2.2}
\end{equation}
\end{lemma}

\bproof Both the statement and its proof are essentially the same as
Lemma 2.2 of \cite{RW}. $~~~\qed$

\section{Counting slopes.}

We first look at some properties of lattices, which we will use to analyze lengths of
short geodesics on flat tori.
A parallelogram $P$ of a lattice $\Gamma$ on $E^2$ is called {\it fundamental} if the vertices of $P$ generate $\Gamma$. 

\begin{lemma} \label{lemma3.1}
Suppose $\Gamma$ is a lattice in the plane
in which the shortest distance between any two vertices is $d\ge 1$ and
the area of a fundamental parallelogram is $\sqrt 3 A^2$, where $A \ge d$. Then there is a
fundamental parallelogram of $\Gamma$ with diameter less than $3A^2$.
\end{lemma} 

\bproof Let $O_1$ and $O_2$ be two independent vertices which have shortest distance to the origin $O$.
For the triangle $OO_1O_2$, let $\alpha$ be the angle at $O$,
and $l_1$, $l_2$ and $l$ be the lengths of $OO_1$, $OO_2$ and $O_1O_2$.
We may assume that $l_2\ge l_1$ and $\displaystyle \alpha \le \frac{\pi}2$, as otherwise we can replace one of the
vertices by its inverse.
Then $l\ge l_2\ge l_1=d$, $\displaystyle \alpha \ge \frac {\pi} {3}$, and therefore
$\displaystyle A^2 \sqrt 3 = l_1 l_2 \sin \alpha \ge l_1 l_2 \frac {\sqrt 3}{2}$, so that $l_1 l_2 \le 2 A^2$.
Since $1 \le d = l_1 \le A \le A^2 $ and $\displaystyle l_2 \le \frac{2A^2}{l_1} \le 2A^2$, 
we have $l_1 + l_2 \le 3A^2$, and the diameter is less than $3A^2$ as claimed.
 $~~~\qed$

Following the notation of Lemma~\ref{lemma3.1},
we have a fundamental parallelogram $P$ of $\Gamma$ spanned by $OO_1$ and $OO_2$ in the Euclidean plane $E^2$.
Letting $h$ be the height of $O_2OO_1$ over $OO_1$, then $dh=A^2 \sqrt 3 \ge d^2 \sqrt 3$
and we obtain:

\begin{lemma} \label{lemma3.2}
 $h\ge d \sqrt 3$.
 \end{lemma}

Next, applying a Euclidean isometry,
we can assume that $e'_1=OO_1=d(1,0)$ and $e'_2=OO_2=d(x,y)$, where 

\begin{equation}
0\le x \le 1/2,\, x^2+y^2\ge 1.\label{eq3.1}
\end{equation} 

By Lemma~\ref{lemma3.2}, we have
\begin{equation}
y = \frac hd \ge \sqrt 3 .\label{eq3.2}
\end{equation}

Now $||ae'_1+be'_2||=d^2((a+bx)^2+by^2)$. So the number of primitive vertices of $\Gamma$ in $D(2g\pi)$, is the cardinality of the set
$$
\{(a,b)\in Z^2: \mbox{gcd} (a,b)=1,\, (a+bx)^2+(by)^2\le (2g\pi/d)^2 \}. $$
Since we are going to find an upper bound,
by condition (\ref{eq3.2}) we may assume that $y=\sqrt 3 $, and then by (\ref{eq3.1}), we have the following 

\begin{lemma} \label{lemma3.3}
$$N(g,d)=\max \{\frac 12 N(g,d,x), 0\le x\le 1/2\},$$ where
$$
N(g,d,x)= \#\{(a,b)\in Z^2: \mbox{gcd} (a,b)=1,\, (a+bx)^2+3b^2\le (2g\pi/d)^2\}.
$$
Moreover $N(g,d,x)$ is a locally constant
function of $x$, with any change of its value
occurring only on the finite set
$E_{a,b}=\{x: (a+bx)^2+3{b^2}= (2g\pi/d)^2\}$.
\end{lemma}

{\bf Remark on Lemma~\ref{lemma3.3}.} Lemma~\ref{lemma3.3} is related to Lemma 12 of \cite{BH}.
We have generalized the calculation of $N(1,1)$ in \cite{BH} to our $N(g,d)$,
and improved the estimate $y\ge \sqrt 3/2$ in \cite{BH},
to $y\ge \sqrt 3$.
The function of two variables $(p,x)$ in \cite{BH} then becomes a function of just one variable
$x$, and the computation is significantly simplified. 

Let $E^2_+$ be the subset of $E^2$ with $y>0$ and $D_+(2g\pi)=D(2g\pi)\cap E^2_+$.
To compute $N(g,d)$, we bound the number of primitive vertices in $D_+(2g\pi)$ on the relevant
lattices, then add one (since there is one primitive vertex lying on the positive $x$-axis).
Let $L_x$ be the lattice generated by
$e_1=(0,1)$ and $e_2(x)=(x, \sqrt 3)$, $0\le x \le 1/2$. 

When changing $L_0$ to $L_{1/2}$ via $L_x$, we see that there are four lattice
points crossing the upper half-circle $\partial D_+(2\pi)$ to enter $D_+(2\pi)$ and
four lattice points crossing $\partial D_+(2\pi)$ to leave $D_+(2\pi)$.
The value of $N(g,d,x)$ takes only three values, 24, 23, 22,
when $0\le x\le 1/2$ and therefore $N(1,1)=24$.
Moreover 24 is reached at $\alpha=\pi/2$ ($x=0$) and $\alpha=\pi/3$ ($x=1/2)$.
If we replace $D_+(2\pi)$ by $D_+(2\pi/1.15)$, we get $N(1,1.15)=18$. 

We list some values of $N(g,d)$ which were obtained by computer calculations: 

\begin{lemma} \label{lemma3.4}
\begin{eqnarray*}
N(1,1)=24,\ N(2,1)=92,\ N(3,1)=198,\ N(4,1)=355,\ N(5,1)=549, \\
N(6,1)=792,\ N(7,1)=1076,\ N(8,1)=1396,\ N(9,1)=1776,\\
N(10,1)=2186,\ N(20,1)=8715,\ 19599\le N(30,1)\le 19600, \\
N(1,1.15094)=18,\ N(2, 1.15094)= 69,\ N(4,1.15094)= 263.
\end{eqnarray*}
\end{lemma}

Finally, we give an asymptotic value for $N(g,d)$. 

\begin{theorem}\label{3.5}
$\displaystyle \lim_{g \to \infty}\frac {N(g,d)} {4\sqrt {3} (g+0.5)^2\pi}\le 1.$
\end{theorem}

\bproof Let $\Gamma$ be a lattice satisfying the conditions of Lemma~\ref{lemma3.1}. By Lemma~\ref{lemma3.1}, we have
a fundamental parallelogram $P$ of diameter $\le 3A^2 $,
where $\sqrt 3 A^2$ is the area of the parallelogram. Then $\Gamma(P)$ gives a tessellation of $E^2$. If a vertex
$V=\gamma(O)$ lies in $D(2g\pi)$ for some $\gamma\in \Gamma$, then $\gamma(P)$ lies in
$D(2g\pi+3A^2)\subset D(2(g+0.5A^2)\pi)$.
Since the area of $D(2(g+0.5A^2)\pi)$ is $(2(g+0.5A^2)\pi)^2\pi$, it follows that there are at most
$\displaystyle \frac {(2(g+0.5A^2)\pi)^2\pi}{\sqrt 3 A^2}$ vertices in $D(2g\pi)$. Since $A$ is a constant and $A \ge 1$, when $g$ is large enough, we have
$$
\frac {(g+0.5A^2)^2}{A^2}\le {(g+0.5)^2}.
$$

By a classical formula due to Dirichlet (pp. 63-64, {\cite{Ap})

$$\lim_{n\to \infty}
\frac {\{\# (r,s) |\mbox{where $r$ and $s$ are coprime, $r^2+s^2\le n^2$} \} } {\{\# (l,m)
|\mbox{where $l$ and $m$ are integers,
$l^2+m^2\le n^2$}\}} =\frac { 6}{\pi^2}. $$
So we have
$$\lim_{n\to \infty} N(g,d)\le \lim_{n\to \infty} \frac 12 \frac {(2(g+0.5)\pi)^2\pi} {\sqrt 3}\frac 6 {\pi^2}
\le \lim_{n\to \infty} 4\sqrt {3}{(g+0.5)^2}\pi.$$ $~~~\qed$

\section{Finiteness for hyperbolic manifolds.}

In this section we apply our calculations to get bounds on the slopes of essential surfaces.

\begin{theorem} \label{4.1}
Suppose $M$ is a compact orientable 3-manifold
with $\partial M$ a torus and that int($M$) admits a complete hyperbolic
metric of finite volume.
Given $g \ge 0$, the number of boundary slopes of an essential immersed surface of
genus at most $g$ is bounded by the function $N(1,d)$ if $g \le 1$ and $N(g,d)+1$ if $g>1$,
where $d\ge 1$ is the cusp length. 
\end{theorem}

Before proving Theorem~\ref{4.1} we discuss some consequences.
A famous result of Thurston, (see \cite{T}, \cite{Go}) shows that if a knot complement contains no essential spheres and tori, then 
at most finitely many 3-manifolds obtained by Dehn surgery on the knot contain essential spheres or tori.
This can be generalized to immersed surfaces of any genus.

\begin{corollary}
Suppose $M$ is an orientable 3-manifold
with $\partial M$ a torus. If $M$ contains no closed
$\pi_1$-injective surfaces of genus at most $g$, then at most $N(g)+1$ surgeries on
the cusp give closed 3-manifolds containing closed $\pi_1$-injective
surfaces with genus at most $g$.
\end{corollary}

We will also prove the following:

\begin{theorem} \label{4.3}
Suppose that $M$ is a compact 3-manifold whose interior admits a complete hyperbolic metric of finite volume.
Then there is a finite collection of boundary slopes $B_i$ for the $i^{th}$ component of $\partial M$
so that if $F$ is any essential immersed surface with genus at most $g$, then one of the
boundary curves of $F$ is contained in some $B_i$.
\end{theorem}

The proofs of Theorem~\ref{4.1} and~\ref{4.3} are based on some results in minimal surface theory,
which extend some standard results in the closed case.
See \cite{HS} for an introduction to the techniques of least area surfaces in 3-manifolds.
We need to use least area existence results in the category of non-compact surfaces.
Such a result is a fairly simple extension of existence results of Schoen-Yau,
but we do not know of an explicit statement in the literature.
We say that a surface is {\it least area}
in its homotopy class if any compact subsurface
is least area in the homotopy class of the subsurface (rel boundary).

\begin{theorem} \label{4.4}
Let $M$ be a compact 3-manifold whose interior admits
a complete hyperbolic metric of finite volume. Let $F$ be an
essential surface in $M$ with finite genus and finitely many
boundary components.
Then int($F$) is properly homotopic in int($M$) to a surface $F'$
which has least area in its homotopy class.\end{theorem}

\bproof Take an exhausting sequence of submanifolds $M_i$ of int($M$),
each of whose boundaries consist of horotori,
so that $M_i \subset \mbox{int}(M_{i+1})$ and $\cup M_i = M$.
We now consider new Riemannian manifolds $M_i'$ obtained from $M_i$ by
altering the hyperbolic metric on $M_i$ so that 

\begin{enumerate}
\item the metric is unchanged on $M_{i-1}$. 

\item A collar of the boundary of $M_i'$ has a flat product metric, with the boundary of $M_i'$ a flat torus. 

\item The sectional curvature of $M_i'$ is non-positive. 
\end{enumerate}

For the construction of such a metric see \cite{Le}.
By a proper homotopy of $F$, it can be arranged that $F$ is transverse
to each $\partial M_i$ and also the intersection of $F$ with each product
region $M_i - \mbox{int} M_{i-1}$ consists of a collection of essential annuli.
We can construct a sequence of smooth closed Riemannian manifolds $N_i'$
by doubling $M_i'$ along its tori boundary components.
The manifold $N_i'$ contains a surface $G_i$ corresponding to the double of $F$.
Since $F$ is essential, Lemma~\ref{lemma2.1} implies that $F$ is injective
on relative $\pi_1$, and its double $G_i$ is a closed essential surface.
The existence result of Schoen-Yau \cite{SY} applies to establish the existence of a least area surface $G_i'$ homotopic to $G_i$.
Theorem 6.7 of \cite{FHS} implies that the intersection of $G_i'$ with the least area torus
$\partial M_i'$ consists of one curve for each boundary component of $F$.
Define $F_i$ to be $G_i' \cap M_i'$.
Then $F_i$ is a least area surface (rel boundary) which is properly homotopic to $F\cap M_i$ in $M_i$.

Now choose a ball $B$ in $M_i$ and consider the intersection of $B$ with the minimal surfaces $F_j$, $j>i$.
The Gauss-Bonnet Theorem for a closed immersed surface $H$ states that
\begin{equation}
\int_H K dA = 2\pi(2-2g) \label{eq5}
\end{equation}
where $g$ is the genus of $H$ and $K$ is the induced curvature on the surface $H$.

Since $G_j' $ is minimal in $N_j$ its induced curvature $K$ is at least as negative as
the sectional curvature of $N_j'$.
Moreover $K$ is smaller or equal to -1 at points where $N_j'$ is hyperbolic.
This implies that the areas of the
intersections $G_j' \cap B$ are uniformly bounded. 
The norm of the second fundamental form, or equivalently the principal curvatures of $G_j' \cap B$
are also uniformly bounded by a result of Schoen \cite{Sc}, since the surfaces
$G_j'$ are least area, and therefore stable.

Given a sequence of least area minimal surfaces with uniformly bounded area and principal curvatures, the
Ascoli-Arzela theorem implies that a subsequence converges in $B$.
It is a property of least area surfaces that such a
limit is also a smooth minimal immersion (see \cite{HS}).
The convergence may be with multiplicity, in which case the
original surface can be recovered by taking appropriate covers of the limit surface. The limit
surface can be extended to all of int($M$) by covering
int($M$) by balls and taking a diagonal subsequence, as in \cite{HS}. This provides the desired surface $F'$.
 $~~~\qed$

{\bf Remark.} We could apply a weaker result than Theorem~\ref{4.4} for our applications.
It suffices to know the existence of a minimal $F_i$ for large enough $i$, allowing us to avoid the
convergence arguments giving a complete minimal surface. A result related to Theorem~\ref{4.4} can be
found in \cite{Ru}.

{\bf Proof of Theorem~\ref{4.3}.}
Choose a collection of maximal horoballs in int($M$), one for each cusp,
so that these horoballs are as big as possible subject to having non-overlapping interiors.
Now push each boundary horotorus slightly outwards towards the cusp it bounds, so that the horotori become
disjoint. Call the resulting horotori $T_i$.
Note that the choice of maximal horoballs is not unique if there is more than one cusp,
but every choice gives a cusp length of at least one in each horotorus.

Next cut off the cusps along each $T_i$ to give a compact, non-complete, hyperbolic
3-manifold $M'$ with boundary a union of flat tori. Clearly $M'$ is homeomorphic to $M$.

Now suppose $F$ is an essential immersed compact surface in $M$ of genus $g$ having $n$ boundary components $c_1,\dots ,c_n$.
Then $F$ is $\partial$-incompressible, by Lemma~\ref{lemma2.1}. Using Theorem~\ref{4.4}, we can properly homotop int($F$) to a least area surface
in int($M$). 
We abuse notation somewhat by also referring to the complete minimal surface we obtain as $F$.

Since $F$ is a complete minimal surface with $n$ cusps, $K<-1$.
The Gauss-Bonnet Theorem gives that
$$
Area(F) (-1) \ge 2\pi \chi(F) = 2\pi (2-2g-n).
$$
So
$$
Area(F) \le 2\pi (2g-2+n).
$$

Fix any positive real number $e$.
We define the collection of boundary slopes $B_i$ so that a boundary slope is in $B_i$ if the
corresponding geodesic on the horotorus in $T_i$ has length less than $2\pi + e$.

We can estimate the area of a surface $F$ in a cusp by using the co-area formula.
This formula implies that the area of $F$ is greater than the integral of the lengths of the intersection
of $F$ with the horotori in the cusp.
For $F$ a boundary incompressible surface which intersects a horotorus $T$
on the boundary of a cusp in geodesics of total length $L$,
the area of $F$ in the cusp is greater or equal to
$$
\int_1^\infty L/y^2 dy = L.
$$

Adding the contribution of the cusp corresponding to each $c_i$,
\begin{equation}
 \Sigma_{i=1}^n L (c_i) \le Area(F) \le 2\pi (2g-2+n). \label{eq4.2}
\end{equation} 

If the boundary of $F \cap T$ contains
$n$ curves all having length greater than $(2\pi + e)$, then $Area(F) > (2\pi + e)n$.
Combining with the Gauss-Bonnet estimate,
we get a contradiction if $e$ is sufficiently large. In particular we cannot have
$e > (g-1)4\pi /n.$ Therefore one of the boundary components must be shorter than $2\pi + e$
and therefore lies in $B_i$. $~~~\qed$

{\bf Remarks.} (1) It suffices to assume that int($M$) has a complete
Riemannian metric with curvature 
less than some constant $C<0$.
However we then need to make some additional arguments to establish existence of a minimal surface.

(2) It suffices to assume that $F$ injects on simple loops and arcs,
as this is all that is needed for the Schoen-Yau existence theorem for least area surfaces. 

From (\ref{eq4.2}) we have
\begin{equation}
\Sigma_{i=1}^n (L (c_i)-2\pi) \le 2\pi(2g-2).\label{eq4.3}
\end{equation} 

The case of a hyperbolic knot complement is of special interest. 

{\bf Proof of Theorem~\ref{4.1}.} In the proof of this and the next pair of theorems,
we assume that our genus $g$ essential surfaces $F$ have been properly homotoped to
least area immersions in the complete hyperbolic metric on int($M$), as in the proof of Theorem~\ref{4.3}. 

Let $c$ be the boundary slope of $F$. Then (\ref{eq4.2}) can be written as 
\begin{equation}
L(c)\le \frac{2\pi(2g-2+n)}n
\label{eq4.4}
\end{equation} 
where $n$ is the number of boundary components of $F$.

If $g=0$, then $n\ge 3$ and we have $L(c)\le 2\pi$. 

If $g=1$ we have $L(c)\le 2\pi$.

If $n=1$, then $c$ is homologically
zero, and there is at most one such slope in $\partial M$.
We have $L(c)\le 4g\pi - 2\pi$ in this case.
Below we assume that $n\ge 2$.

If $g>1$, since $n\ge 2$ we also have $L(c)\le 2\pi(g-1)+2\pi=2g\pi$. 

In conclusion, with the possible exception of the null-homologous slope.
$L(c)\le 2g\pi$ when $g>1$. 

Since the horotorus where $c$ lies can be arbitrarily close to the maximum horotorus,
we assume for convenience that $c$ is actually contained in the maximum horotorus.
The universal cover of the maximum horotorus is $E^2$, which is tessellated
by the fundamental domain of the maximum horotorus, namely a parallelogram $P$.
We may assume that a vertex of $P$ is at the origin.
By an observation of Colin Adams \cite{Ad1}, such a $P$ contains two disjoint discs of
diameter $d$, neighborhoods of the two points where the maximal
cusp first touches itself, and its area is at least $\displaystyle d^2\sqrt 3$. 
The number of boundary slopes of length at most $2g\pi$ is half of the number
of primitive vertices of the tessellation lying in $D(2g\pi)$,
since two vertices $(p,q)$ and $(-p,-q)$ correspond to the same slope.
This is bounded above by $N(g,d)$, for $d$ the smallest distance between two vertices in the
tessellation of $E^2$.

It follows that the number of simple closed 
geodesics on the horotorus with length at most $2g\pi$ is bounded by $N(g,d)$.
All slopes except the single slope arising from the case $n=1$ are covered by 
this case, so the total number of slopes is at most $N(g,d)$+1 when $g>1$, and 
Theorem~\ref{4.1} is proved. $~~~\qed$

\begin{theorem} \label{4.5}
Let M be a compact 3-manifold with boundary a torus
whose interior admits a complete hyperbolic metric of finite volume. 
Suppose $\alpha_i$, $i=1,2$, are two different slopes on $\partial M$ such that $\alpha_i$ is a boundary slope of a punctured $\pi_1$-injective immersed
surface $F_i$ of genus at most $g_i>0$.
If each $\alpha_i$, $i=1,2$ is not homologically zero in $M$,
then the intersection number
$\Delta (\alpha_1,\alpha_2)$
is bounded by $11.8g_1g_2$. If one of the
$\alpha_i$ is homologous to zero, then the bounds
above should be doubled.\end{theorem}

\bproof
First note that at most one slope of a simple closed curve on $\partial M$ is
homologous to zero in $M$.

Suppose each $\alpha_i$, $i=1,2$, is not homologous to zero.
Then $F_i$ has at least two boundary components, and by the
proof of Theorem~\ref{4.1}, the length
of $\alpha_i$ in the maximal cusp is at most $2g_i \pi$.
Hence the area of the parallelogram $P$ spanned by the lifts of $\alpha_1$
and $\alpha_2$ on the Euclidean plane at height 1 is bounded by $4\pi^2 g_1g_2$. Let $T$ be the area of the boundary of the maximal cusp. The intersection number $\Delta (\alpha_1, \alpha_2)$ is the area of $P$ divided by the area of $T$, that is 

\begin{equation}
\Delta (\alpha_1, \alpha_2)=\frac{\mbox {area $P$}}{\mbox{area $T$}} \le
\frac{4\pi^2 g_1 g_2}{\mbox {area $T$}}.\label{eq4.5}
\end{equation} 

The area of $T$ is at least $3.35$ by Lemma 2.2.
We calculate $\displaystyle \frac{4\pi^2}{3.35}\le 11.8$.
It follows that $\Delta (\alpha_1,\alpha_2)$
is bounded by $11.8g_1g_2$.

If some $\alpha_i$, say $\alpha_1$, is homologous to zero,
then by the proof of Theorem~\ref{4.1},
the length of $\alpha_i$ in the maximal cusp is
at most $4g_i \pi$ and $\alpha_2$ is not homologous to zero.
So the above bounds should be doubled to bound the intersection number.
This proves the Theorem.
 $~~~\qed$

{\bf Remarks.}
We note that without the area estimate due to Cao and Meyerhoff we
would have the weaker estimate of $17g_1g_2$ instead of
$11.8g_1g_2$. We also note that
Agol \cite{Ag} has recently obtained estimates on the length
of slopes of punctured spheres and tori which give constants sharper than the ones obtained here for those cases.
Agol's methods can be combined with ours to show that the constant
$2\pi$ can be improved to six in formulas such as
(\ref{eq4.2}). See also Lackenby \cite{La}.
If we apply these in Theorem 4.5, $11.8g_1g_2$ can be
replaced by $10.8g_1g_2$, and in Lemma 3.4
we have
\begin{eqnarray*}
N(1,1)=22,\ N(2,1)=80,\ N(3,1)=182,\ N(4,1)=323,\ N(5,1)=502, \\
N(6,1)=721,\ N(7,1)=979,\ N(8,1)=1277,\ N(9,1)=1616,\\
N(10,1)=1994,\ N(20,1)=7955,\ 17874\le N(30,1)\le 17875 , \\
N(1,1.15094)=16,\ N(2, 1.15094)= 62,\ N(4,1.15094)= 242.
\end{eqnarray*}

We next show how to get a bound for the number of boundary slopes.

\begin{theorem} \label{4.6}
Suppose $M$ is a compact orientable 3-manifold
whose interior admits a complete hyperbolic metric of finite volume. 
Suppose $F\to M$ is an essential immersion of a surface of genus $g$ with boundary slope $\lambda$. Then
$\#\partial F \le (g-1)C$ for some $C>0$, with at most 24 exceptions on $\lambda$.
Given any $k>1$, then $\displaystyle \#\partial F\le \frac {2g-2}{k-1}$, with at most $N(k,1)$ exceptions.
\end{theorem}

\bproof We saw in Lemmas~\ref{lemma3.3} and \ref{lemma3.4} that
all primitive lattice points in $D(2g\pi)$ lie in $D(2g\pi)-D(2\pi)$,
with at most 48 exceptions, have distance $\ge 2\pi(1+e)$ from the origin,
where $e$ is a positive number, independent of $M$.
Let $C=2/e$. By applying (\ref{eq4.3}), we have
$$
2n\pi e\le 2\pi (2g-2),
$$
i.e. $n\le C(g-1)$.

Similarly, all primitive lattice points in $ D(2g\pi)$ are contained in $D(2g\pi)-D(2k\pi)$
with at most $N(k,1)$ exceptions.
So they have distance $\ge 2k\pi$ from the origin,
and we get
$$n(2(k-1)\pi)\le 2\pi(2g-2)$$ so that
$\displaystyle n\le \frac {2g-2}{k-1}.$ $~~~\qed$

Suppose $K$ is a hyperbolic knot in $S^3$.
Let $(K,\lambda)$, $c(K, \lambda, g)$ and $n(K, \lambda, g)$ be defined as in the introduction. 

\begin{theorem} \label{4.7}
If $(K,\lambda)$ is irreducible, then
$c(K, \lambda, g)\le n(K,\lambda,g).$
\end{theorem}

\bproof Suppose $F$ is an essential surface of genus $g$ in $M = (K, \lambda)$,
the manifold obtained by $\lambda$-surgery on $K$,
and that $F$ has been homotoped so that it intersects the core of the surgery transversely and
realizes the minimal geometric intersection number $m$ with the core. Then $m \le c( K,\lambda, g)$.
Moreover we can assume that $m = c(K, \lambda, g)$ if $c(K,\lambda,g)$
is finite. Note that $m$ can be arbitrarily large if $c(K,\lambda,g)$ is infinite,
for an appropriately chosen $F$.

We make $F$ least area in its homotopy class in a metric constructed
by blowing up the metric in a small regular neighborhood of $K$.
Then the least area map homotopic to $F$ does not increase the minimal intersection number with $K$. 
By \cite{FHS} if $F$ is least area, it lifts to an embedding $\tilde F$ in the covering space $M_F$ of $(K, \lambda)$
with fundamental group given by the subgroup $\pi_1(F)$ in $\pi_1(M)$.

Remove a regular neighborhood int$(N(K))$ of the core and denote $F- \mbox{int} (F \cap N(k))$ by $F_1$.
$F_1$ is a proper immersed surface of
genus $g$ with $m$ boundary components.
We will show that $F_1$ is essential, so that $n(K,\lambda, g) \ge m$ which in turn implies that $n(K,\lambda,g) \ge c(K,\lambda,g)$. 

Next remove all the lifts of $\mbox{int}(N(K))$ from $M_F$. This gives a
covering space of $M-\mbox{int} (N(K))$. Notice that $F_1$ lifts to an embedding
$\tilde F_1$ in the covering space.
If $F_1$ is not essential in $M-\mbox{int}(N(K))$, then $\tilde F_1$ is not essential and
by the loop theorem,
there is an embedded compressing disc for $\tilde F_1$. The boundary of this disk
projects to give an essential simple loop
$c$ on $F_1$ which bounds a singular disk in $M-\mbox{int}(N(K))$.
The curve $c$ bounds a singular disc on $F$, which must meet $K$. Since $(K, \lambda )$ is irreducible,
the union of these two immersed discs represents a null-homotopic 2-sphere,
and we can homotop $F$ in $(K,\lambda )$ to move the disk on $F$ bounded by $c$ to the singular disk in
$M-\mbox{int} (N(K))$, and thus we reduce the number of intersections between $K$ and $F$, giving a contradiction. $~~~\qed$

{\bf Remark.} For any hyperbolic knot in $S^3$, it is known that there is at most one reducible surgery,
and the cabling conjecture states that there is no reducible surgery.

By Theorem~\ref{4.6}, Theorem~\ref{4.7} and the remark, we have 

\begin{corollary}
For any hyperbolic knot $K$ in $S^3$ and any $g>1$, $c(K, \lambda ,g) \le (g-1)C$ for some $C>0$, with at most 25 exceptions for $\lambda$.
Also given $k>1$, $\displaystyle c(K, \lambda ,g) \le \frac {2g-2}{k-1}$, with at most $N(k,1)+1$
exceptions for $\lambda$.
\end{corollary}

\section{Finiteness for Haken manifolds.} 

In this section we discuss the case of a general (possibly
toroidal) Haken manifold with boundary.

\begin{theorem} \label{5.1}
Suppose $M$ is an orientable Haken 3-manifold with $\partial M$ a torus.
Then there are only finitely many boundary slopes realized by orientable
essential proper surfaces of genus at most $g$.
\end{theorem}

\bproof
Let $\Gamma$ be the Jaco-Shalen-Johannson decomposition tori of $M$. If
$\Gamma=\emptyset$, then $M$ is either a hyperbolic
3-manifold or a Seifert manifold.
If $M$ is a hyperbolic manifold, then the conclusion of Theorem~\ref{5.1}
follows from Theorem~\ref{4.1}.
If $M$ is a Seifert manifold, then the boundary slope is unique by (\ref{eq2.2})
in Lemma~\ref{lemma2.3}.

Below we assume that $\Gamma $ is not empty. Call each component of
$\overline{M-N(\Gamma)}$ a {\it vertex manifold}, where $N(\Gamma)$ is a
regular neighborhood of $\Gamma$.

Let $M_*$ be the vertex manifold containing the boundary torus of $M$.
Suppose there are infinitely many boundary slopes $\{B_n\}$ for essential
immersed surfaces of genus at most $g$. Then for each $B_n$, there is an
essential surface $F_n$ of genus at most $g$ such that $\partial
F_n$ has $l_n$ components, each with slope $B_n$.
First deform $F_n$ so that the number of components of $F_n\cap\partial
N(\Gamma)$ is a minimum. Let $F^*_n$ be the union of the components of
$F_n\cap
M_*$ with boundary components on $\partial M$.

Let $l^*_n$ be the number of boundary components of $\partial F^*_n$ on
${\partial M_*}-\partial M$.

\begin{lemma} \label{lemma5.2}
For any constant $C>0$, there is a constant $C'$ so that $l^*_n< Cnl_n$, 
whenever $n > C'$.
\end{lemma}

\bproof
Let $S_n$ denote all the components of $F_n- F^*_n$. Since both
$M_*$
and $M-M_*$ are boundary irreducible, no component of $S_n$ is a disc.
To recover $F_n$ from $F^*_n$ and $S_n$, we identify
the loops of $\partial S_n$ and $\partial F^*_n$ in three steps:

\begin{enumerate}
\item Identify a minimum number of pairs of loops of $\partial S_n$
and $ \partial F^*_n$ to
form a connected surface denoted $F_n'$, which contains all the boundary 
curves of $F_n$. Let $S_n'=F_n-F_n'$.

\item Glue each
component of $S'_n$
which has more than one boundary component to $F_n'$ along
exactly one of its boundary curves to form
$F_n''$ and let $S_n''=F_n-F_n''$;

\item Identify
all the remaining pairs of loops of $\partial S_n''$ and
$\partial F''_n$ to get back $F_n$.
\end{enumerate}

Assume that for some constant $C>0$, $l^*_n > Cnl_n$ for some
unbounded sequence of choices for $n$, so that no
constant $C'$ exists as required. We can suppose that
$n$ is chosen arbitrarily large in this sequence and seek a
contradiction.
Since $F^*_n$ has at most $l_n$ components, the maximum
number of components of $\partial S_n$ and $\partial F^*_n$ which
are identified is $2l_n$. So there are at least $(Cn-2)l_n$ components
of $\partial S'_n$. In the second step the number
of boundary components of $S'_n$ we glued is no more than $1/2(Cn-2)l_n$,
and
therefore $\partial S''_n$ has at least $1/2(Cn-2)l_n$ components.
So in the third step we need to identify at least $1/2(Cn-2)l_n$ pairs of
components of $\partial S_n''$ and $\partial F''_n$. The surface
genus increases by one
when we identify such a pair, so the genus of $F_n$ is at least
$\displaystyle 1/2(Cn-2)l_n$,
which is unbounded, since $n$ can be chosen arbitrarily large.
So this contradiction proves the lemma. $~~~\qed$

Since the genus of $F_n$ is assumed to be at most $g$,
the genus of $F^*_n$ is also at most $g$.
By Lemma~\ref{lemma5.2}, to prove Theorem~\ref{5.1},
we find an (unbounded) sequence of values of $n$, so that
$l^*_n> Cnl_n$ for some constant $C>0$.

We have two cases.

Case (1) $M_*$ is hyperbolic. Up to a choice of subsequence
of $n$, we may assume that
the length of $B_n$ is larger then $2(n+1)\pi$. Moreover by 
(\ref{eq4.2}), we have

$$ 2(n+1)l_n\pi + \sum_{c\in \partial M_*-\partial M} L(c)\le$$

$$\sum_{c'\in \partial M} L(c') + \sum_{c\in \partial M_*-\partial M}L(c)
\le
2\pi(2g-2+l^*_n+l_n).$$

So
$$nl_n\le 2g-2+l^*_n,$$
i.e., $l^*_n> Cl_nn$, for $C=2$ and $n>2g-2$.

Case (2) $M_*$ is a Seifert manifold and $B_n=(u_n,v_n)$. Let $O(M_*)$ be
the Seifert orbifold for $M_*$. Denote the Euler characteristic of
$O(M_*)$ by $\chi_*$.
There are two subcases.

Case (2a) $u_n$ is unbounded as $n \rightarrow \infty$. We may assume that
$u_n>n$, by taking a subsequence of values of $n$.

Notice that the projection $p: F^*_n\to O(M_*)$ is an orbifold branched
covering
of degree at least $l_n u_n$. In fact, by Hass \cite {Ha} we may assume
that $F^*_n$ is horizontal relative to the Seifert fibering. Recall also
that $F_n$ has $l_n$ boundary curves, each of which has coordinates
a non zero multiple of $(u_n, v_n)$.
Using the estimate of the degree of $p$, it follows that
$$l_n u_n \chi_* \ge \chi(F^*_n)=2\# F^*_n-2g(F^*_n)-\#\partial F^*_n.$$

So

\begin{equation}
 \# \partial F^*_n \ge
-l_n n\chi_*+2\# F^*_n-2g(F^*_n) \ge
-l_n n\chi_*-2g.\label{eq5.1}
\end{equation}
By \ref{5.1}, there are at least
$l_n(-n{\chi_*}-1)-2g$ components of $\partial F^*_n$ on $\partial
M_*-\partial M$. Now since $\chi_* <0$, as $n$ tends to infinity,
we see that $l^*_n>Cnl_n$, where $C = 1/2$.

Case (2b) $|u_n |$ is bounded by a constant $u>0$, so $|v_n|$ tends to
infinity with increasing $n$. We may assume that $|v_n|>2n$ by
choosing a subsequence of values of $n$.

For convenience, the coordinates of a closed curve
$c\subset T(\mu,\lambda )$ will be denoted by $(u_c, v_c)$. Let $C_{n,j}$
be all components of $\partial F^*_n$ lying in $T_j$, where $\partial
M^*=\{T_1,\dots , T_h\}$ and $\partial M=T_h$.

By (\ref{eq2.2}), we have
$$\sum_{j=1}^{h-1}\sum_{c\in C_{n,j}}
\frac {v_c}{u_c} =-l_n\frac {v_n}{u_n}.$$
It follows that
$$
\sum_{j=1}^{h-1}\sum_{c\in C_{n,j}} |v_c|\ge l_n\frac {|v_n|}{u}.
$$
So there is at least one $j$, say $j=1$, such that 
\begin{equation}
\sum_{c\in C_{n,1}} |v_c|\ge
\frac{l_n|v_n|}{(h-1)u}.
\label{eq5.2}
\end{equation}

Since each component of $F^*_n $ is a $\pi_1$-injective surface in the
Seifert fibered
manifold $M_*$ and is not vertical, it must be horizontal. By (\ref{eq2.1}), there
are at most $l_nu$ components of $\partial F^*_n$ lying in $T_1$.
So by (\ref{eq5.2})
the average of $|v_c|$ for curves of $\partial F^*_n$ on $T_1$ is at
least $\displaystyle \frac
{|v_n|}{(h-1)u^2}$.

Let $ C_n$ be the collection of
components of $\partial F^*_n$ lying in $T_1$ such that $\displaystyle
|v_c|\ge\frac {n}{(h-1)u^2}$ .

There are at most $l_nu$ components in $\partial F^*_n \cap T_1$ and so
also in $C_{n,1}-C_n$.
The value of $|v_c|$ for each component in
$C_{n,1}-C_n$ is at most $\displaystyle \frac {n}{(h-1)u^2}$ and $|v_n|>2n$ so
by (\ref{eq5.2}) we have

$$
\sum_{c\in C_n }|v_c| + \frac {nl_n}{(h-1)u}> \sum_{c\in C_n
}|v_c|+\sum_{c\in
C_{n,1}-C_n }|v_c|
>\frac {l_n|v_n|}{(h-1)u}>
\frac {2nl_n}{(h-1)u}.$$
That is
\begin{equation}
\sum_{c\in C_n }|v_c| > \frac {nl_n}{(h-1)|u|}> Cl_nn,
\label{eq5.3}
\end{equation}
where
$\displaystyle C=\frac 1{(h-1)u}$.

Let $M'$ be the vertex manifold of $M$ sharing the torus $T_1$ with $M_*$.
Denote the copy of $T_1$ on $M'$ by $T'_1$ the gluing map by $g:T_1\to
T'_1$.

If $M'$ is hyperbolic,
we assume that

(1) the hyperbolic structure $M'$ is obtained by removing a maximal torus
cusp
from its unique complete finite volume hyperbolic structure,

(2) a Euclidean coordinate system is chosen on $T'_1$,

(3) $g$ is affine.

By (\ref{eq5.3}), when $n$ is
sufficiently large, there is a constant $C$ such that on $T'_1$ we have

\begin{equation}
\sum_{c\in C_n} L(g(c))>Cn l_n. \label{eq5.4}
\end{equation}

If $M'$ is a Seifert manifold, the gluing map $g: T_1(\mu, \lambda )\to
T_1'(\mu', \lambda ')$ is determined by a 2 by 2 matrix
$\displaystyle
A = \left( 
\begin{array}{cc}
p & q \\
r & s
\end{array}
\right)
$,
where $r\ne 0$, and $qr-ps=1$. Let $g_*$ be the induced map on
homology, so that $g_*(u\mu + v \lambda )=u'\mu' + v' \lambda'$. Then $u'=pu+rv$ and
$v'=qu+sv$. Note $|r|\ne 0$, $u_c$ is bounded and $|v_c|>2n$. By (\ref{eq5.3}), when $n$
is sufficiently large, there is a constant $C$ such that on $T'_1$
\begin{equation}
\sum_{c\in C_n} |u_{g(c)}|>C l_nn. \label{eq5.5}
\end{equation}

In either case we must have $M_*\ne M'$.

Let $c \in C_n$ and let $F^{**}_n\subset M^*\cup N(T_1)\cup M'$ be a
subsurface of $F_n$. $F^{**}$ is composed
of $F^*_n$, components of $F_n\cap M'$ which have $g(c)$ as
boundary components and those annuli in $N(T_1)$ connecting $c$ and
$g(c)$ for all choices of $c$. Then clearly

(1) $F^{**}_n$ has at most $l_n$ components.
The genus of $F^{**}_n\cap M'$ must be bounded by $g$. By (\ref{eq5.4})
and the calculation in Case (1) when $M'$ is hyperbolic, or by (\ref{eq5.5}) and
the calculation in Case (2a) when $M'$ is Seifert fibered, it follows that
$\# \partial (F^{**}_n\cap M')-\#\{g(c), c\in C_n\}>C''n l_n$ for some
non-zero constant $C''$.
Consequently,

(2) $\partial F^{**}_n-\partial M$
has at least $C''nl_n$ components.

By (1) and (2), we can apply the proof of Lemma~\ref{lemma5.2} to $F^{**}_n$ to get
that
the
genus of $F_n$ is unbounded when $n$ increases. $~~~\qed$

\section{Acknowledgements.}
We would like to thank C. Adams, C. Gordon, J. Luecke and Y. Wu for helpful conversations,
and X. Zhao for writing a program to compute $N(g,d)$. 
An early version of the paper was written when the three authors
were visiting MSRI, Berkeley in 1996-97. Further work was carried out by the second and third
authors at the Morningside Center, Beijing in 1998.
Research at MSRI is supported in part by NSF grant DMS 9022140.
The first author was partially supported by NSF grant DMS-9704286,
the second by an ARC grant
and the third by an Outstanding Youth Fellowship of NSF China grant.

\begin{flushleft}
Joel Hass, Department of Mathematics,
University of California, Davis, CA 95616.
email: hass@math.ucdavis.edu\\
J.Hyam Rubinstein, Department of Mathematics,
The University of Melbourne, Parkville, Victoria 3052, Australia.
email: rubin@ms.unimelb.edu.au\\
Shicheng Wang, Department of Mathematics,
Peking University, Beijing 100871, China. email: swang@sxx0.math.pku.edu.cn \\
\end{flushleft}

\end{document}